\documentclass[12pt,draftcls,onecolumn]{IEEEtran}
\usepackage{cite}

\ifCLASSINFOpdf
   \usepackage[pdftex]{graphicx}
\else
   \usepackage[dvips]{graphicx}
   \DeclareGraphicsExtensions{.eps}
\fi

\usepackage[cmex10]{amsmath}
\usepackage{amssymb}  
\usepackage{bm}
\usepackage{color}
\usepackage{mdwmath}
\usepackage{mdwtab}
\usepackage{url}
\usepackage{enumerate}
\usepackage[linesnumbered,ruled]{algorithm2e}

\newtheorem{theorem}{\textbf{Theorem}}
\newtheorem{definition}{\textbf{Definition}}
\newtheorem{lemma}{\textbf{Lemma}}
\newtheorem{remark}{\textbf{Remark}}
\newtheorem{assumption}{\textbf{Assumption}}

\newtheorem{alg}{\textbf{Algorithm}}

\definecolor{temp}{rgb}{0,0,0}

\begin{document}
%
\title{Distributed Model Predictive Control of Spatially
Interconnected Systems Using Switched Cost Functions}
%
%
%
\author{Peng~Liu,
        and~Umit~Ozguner
\thanks{This work was supported in part by National Science Foundation under Cyber-Physical Systems (CPS) Program (CNS-1446735), and U.S. Department of Transportation under University Transportation Center (UTC) Progarm (DTRT13-GUTC47).}
\thanks{P. Liu and U. Ozguner are with the Department
of Electrical and Computer Engineering, The Ohio State University, Columbus,
OH, 43210 USA (e-mail: liu.3193@osu.edu; ozguner.1@osu.edu).}
\thanks{Color versions of one or more of the figures in this paper are available online
at http://ieeexplore.ieee.org.}}
%
%

\markboth{IEEE~TRANSACTIONS~ON~AUTOMATIC~CONTROL,~JUNE~2016}%
{Liu \MakeLowercase{\textit{et al.}}: Distributed model predictive control of spatially interconnected systems using switched cost functions}
%



\maketitle

\begin{abstract}
This note proposes a distributed model predictive control (DMPC) scheme with switched cost functions for a class of spatially interconnected systems with communication constraints. Non-iterative and parallel communication strategy is considered to ensure that all distributed controllers complete input updates at each single information exchange step. The proposed DMPC scheme switches the optimization index on a switching surface generated by control invariant sets. With the index-switching strategy, stability of the origin is ensured by a terminal control law. Convergence conditions of the optimal cost to zero are established taking into account the causal link between the presumed trajectory and the optimized trajectory of the previous step. The compatibility constraints preserve the quadratic program property that is desired in practical applications. It is also observed that the proposed DMPC scheme has benefits dealing with systems that need to take into account safety-related spatial constraints.
\end{abstract}


%
\IEEEpeerreviewmaketitle


\section{Introduction}
\IEEEPARstart{S}{patially} interconnected systems have wide application areas ranging from automated highway systems, to unmanned aerial vehicle (UAV) formation, to distributed mobile sensing agents\cite{murray2007recent,tang2005motion}. In a spatially interconnected system, cooperation of distributed subsystems (agents) is usually required to achieve global objectives of the overall system such as behavioral consensus, task allocation, economic savings, and so forth. Cooperation of subsystems could be either systemwide or local depending on the communication network and cooperation strategy. While centralized control can deal with systemwide interactions by controlling all subsystems via a single agent, the difficulty of systemwide organization and maintenance via the central agent deters application of centralized control schemes to distributed systems \cite{stewart2010cooperative,camponogara2002distributed}. Research on spatially distributed control schemes that are flexible, scalable, and "plug-and-play" has been attracting increasing attention in cyber-physical system study.

Distributed model predictive control (DMPC) is one of the distributed control schemes developed to regulate spatially interconnected systems subject to system and communication constraints. DMPC inherits properties from model predictive control (MPC) of taking explicit account of state and input constraints. The distributed controllers cooperatively solve a constrained optimal control problem (COCP) in a receding horizon fashion, and implement respectively first items of optimized input sequences to corresponding subsystems. Different from decentralized MPC, DMPC implements an information exchange process via communication connections to fulfill distributed optimization. In general, DMPC schemes can be grouped into different types according to the optimization procedure, process features, communication architecture, and so forth \cite{negenborn2014distributed,maestre2013distributed}. For monolithic systems comprising distributed local controllers, iterative DMPC implementing various distributed optimization approaches has been studied (see, e.g., \cite{venkat2005stability,cai2014rapid,conte2012distributed,stewart2010cooperative,trodden2017distributed}). Similar DMPC schemes employing a coordination layer to allocate distributed optimization tasks have also been reported in \cite{scattolini2009architectures}. For spatially interconnected systems where coupling is usually raised by cost functions and state constraints, non-iterative DMPC strategies have been investigated implementing sequential \cite{richards2007robust,muller2012cooperative}, alternative \cite{trodden2014cooperative,liu2017non-iterative}, and parallel \cite{dunbar2006distributed,dunbar2012distributed} information exchange strategies over neighbor-to-neighbor communication. \textcolor{temp}{In cases where channel reoccurrence rate is not considerably high compared with control rates, e.g. 10 Hz rate of the control channel in DSRC/WAVE\footnote{refer to IEEE 1609.4. It is different from data transmission rate.} for intelligent transportation systems, the distributed control input is expected to update at each data-exchange cycle. Therefore, parallel communication and synchronous timing (c.f. \cite{negenborn2014distributed}) are desirable. Since optimized state trajectories of neighbors are unavailable at the same step to each subsystem implementing parallel information exchange, presumed trajectories of neighbors are used as heuristics in each distributed optimization process. Typically, compatibility constraints \cite{dunbar2006distributed} are introduced for each subsystem to restrain the mismatch between the optimized trajectory and the presumed trajectory. Compatibility constraints are designed together with common MPC conditions to guarantee convergence of the closed-loop trajectory and stability of the equilibrium. Convergence of the closed-loop trajectory to an arbitrarily small level set was investigated in \cite{dunbar2006distributed} by analyzing the overall cost as a function of the overall state. Lyapunov-based approaches can be adopted if further assumption is made on continuity of the control law in state at the equilibrium point\cite{dunbar2012distributed,li2016distributed}. Non-iterative DMPC implementing an input-output form was reported in \cite{liu2016distributed}. }


\textcolor{temp}{As the presumed trajectories can have a causal link to the optimized trajectory of the previous step, analysis of the control law as a function of both the current state and states in the presumed trajectory could provide more insights into non-iterative DMPC schemes implementing compatibility constraints.} This note proposes a non-iterative DMPC scheme using switched cost functions. The idea follows the dual-mode MPC approach and adds an explicit instruction to switch to the terminal control law. By switching the control law and the cost function on a switching surface generated by control invariant sets, stability of the overall system is guaranteed by the terminal control law that does not require compatibility constraints. Convergence of the optimal cost is then addressed by analyzing convergence of the sequence of an augmented state comprising current state, presumed states, and the optimal input. Convergence of the optimal cost indicates  convergence of the closed-loop trajectory to the origin. In addition, decision consensus is achieved via neighbor-to-neighbor communication if the optimal cost converges to zero. Furthermore, the DMPC scheme preserves quadratic program properties by carefully selecting the presumed trajectory and the compatibility constraints.

The rest of this paper is outlined as follows. Section II gives the preliminaries to formulate non-iterative DMPC of spatially interconnected systems. Section III presents the proposed DMPC scheme implementing switched cost functions. Stability conditions are established, and spatial properties are analyzed. Section IV gives a numerical example to illustrate the effectiveness of the proposed DMPC scheme. Section V concludes this work.

\subsection*{Notations}
Throughout this work, let $||\cdot||$ and $||\cdot||_{\infty}$ denote the Euclidian norm and the infinity norm, respectively. $||x||_P^2$ denotes the quadratic form of $x$, i.e. $||x||_P^2=x^TPx$. $\Phi(k,x_0)$ denotes the state trajectory of $x$ at time $k$ with initial condition $x_0$. $\mathcal B_{\rho}(\bm 0)\triangleq \{x\in \mathbb R^n:||x||\leq\rho\}$. $\mathbb Z^+$ denotes the set of nonnegative integers. \textcolor{temp}{$DJ(x)$ denotes the first forward difference of $J$ along $x$, i.e. $DJ(x(t))=J(x(t+1))-J(x(t)), t\in\mathbb Z^+$.} $A^T$ denotes the transpose of matrix $A$. $I$ denotes the identity matrix of appropriate dimention. The superscript $*$ denotes the optimizer or optimal cost according to the context.

\section{Preliminaries}
\subsection{Spatially Decoupled Systems}

\textcolor{temp}{This work focuses on spatially interconnected systems consisting of dynamically decoupled subsystems that have neighbor-to-neighbor communication capability. With the assistance of inter-subsystem communication, the overall system is expected to achieve cooperative properties such as synchronization, formation, string stability, etc., over all subsystems. We focus on the case where cooperation is addressed via coupled cost functions with decoupled input and state constraints. Extensions of the proposed DMPC scheme to coupled state constraints will be discussed in section III-B.} Consider a spatially interconnected system consisting of $M$ identical subsystems with the following discrete-time LTI dynamic
\begin{align}
\label{eqn:all-dyna}
\zeta^{i}(t+1) &=A \zeta^{i}(t)+B \varrho^i(t)
\end{align}
\textcolor{temp}{where $\zeta^i\in \mathbb R^n$, $\varrho^i\in\mathbb R^m$ are the state and the input of subsystem $i$, respectively, with initial condition $\zeta^i(0)=\zeta_0^i$. $(A,B)$ is a controllable pair and $B$ is of full column rank. For systems such as automated highway platooning, UAV formation, the overall system is typically asked to track a desired trajectory. To simplify the analysis of all subsystems in a unified form, a reference frame is introduced. Let $\hat \zeta(t)$ and $\hat \varrho(t)$ denote the reference state and reference input at $t$, respectively. $\hat\zeta(t+1) = A\hat\zeta(t) + B\hat\varrho(t)$. Further define $x^i=\xi^i-\hat\xi^i$, $u^i=\varrho^i-\hat\varrho^i$. The error dynamic is given as follows}
\begin{align}
\label{eqn:dynamic}
x^{i}(t+1) &=A x^{i}(t)+B u^i(t)
\end{align}
with initial condition $x^i(0)=\zeta^i(0)-\hat\zeta^i(0)$. For the rest of this work, we focus on the error dynamics. The state and the input are subject to the following constraints
\begin{align}
x^i(t) & \in\mathcal X^i\subset \mathbb R^n\nonumber\\
u^i(t) & \in \mathcal U^i\subset\mathbb R^m, \quad\forall t\in \mathbb Z^+ \nonumber
\end{align}

The overall system is represented by an undirected graph structure by associating the $i$th subsystem to the $i$th vertex of the graph. Then, the interaction of system $i$ and $j$ is presented by an edge $(i,j)$ in the graph. To better describe information exchanges and interactions between subsystems, we first introduce the following definitions.
\begin{definition}[Undirected graph]
An undirected graph $\mathcal G=(V,E)$ consists of a vertex set $V=\{v_1,\dots, v_M\}$ of $M$ vertices and an edge set $E\subset V\times V$ of unordered pairs $\{e_{ji}=e_{ij}\triangleq (v_i,v_j), v_i,v_j\in V\}$.
\end{definition}

For distributed systems employing asymmetric information exchange structures, e.g. leader-following mobile agents, directed graphs may be applied to formulate the overall system. Without loss of generality, information exchanges between distributed subsystems are assumed to be symmetric throughout this work.
\begin{definition}[Neighbor set]
The neighbor set of a vertex $v_i\in \mathcal G$, is $\mathcal N^i\triangleq \{v_j\in V:(v_i,v_j)\in E, j\neq i\}$.
\end{definition}

The subsystems are interconnected by the means of coupled cost functions
\begin{equation}
\label{eqn:cost-general}
J_i\triangleq J(x^i,u^i,\{x^j,u^j\}),\quad x^j\in \mathcal N^i
\end{equation}

\textcolor{temp}{Throughout this work, the undirected graph $\mathcal G$ of the spatially interconnected system is assumed to have a spanning tree, i.e., there is a path between two arbitrary vertexes in $\mathcal G$.}
\textcolor{temp}{This assumption indicates that decision consensus of all subsystems can be made through multiple information exchanges even if the system is not all-to-all connected.} When $\mathcal N^i=\{x^j:j=1,\dots, M\}$, the distributed control problem degrades to a centralized optimal control problem. In general, the graph edges representing subsystem interactions can be time-varying. For distributed control of spatially decoupled systems, the emphasis is on coordination of subsystems between high-level decision-making stages of the overall system. For the sake of simplicity, the graph structure is assumed to be time invariant throughout this note. To further analyze stability properties of DMPC, the following definitions are presented.

\begin{definition}
$\alpha :[0,a)\to [0,\infty)$ is said to be a class $\mathcal K$ function, if it is continuous on $[0,a)$ and strictly increasing with $\alpha(0)=0$.
\end{definition}

\textcolor{temp}{A class $\mathcal K$ function is further said to belong to class $\mathcal K_{\infty}$ if i) $a = \infty$, and ii) $\lim_{x\to\infty}\alpha(x) = \infty$.}

\begin{definition}[Multi-parametric quadratic programs]
A multi-parametric quadratic program (mp-QP) is a multi-parametric program in the following form
\begin{equation}
\label{eqn:mpqp}
\begin{aligned}
J^*(x) &= \min_{z} J(z,x)=\frac{1}{2}z^THz\\
& s.t. \quad Gz\leq w+Sx
\end{aligned}
\end{equation}
where $z\in \mathbb R^s$ and $x\in\mathbb R^n$ are the optimization variable and the parameter, respectively. $J(z,x):\mathbb R^{s+n}\to \mathbb R$, $H\succ 0$, $G\in\mathbb R^{p\times s}$, $w\in\mathbb R^{p}$, and $S\in \mathbb R^{p\times n}$.
\end{definition}

\subsection{Distributed Model Predictive Control}
A non-iterative DMPC scheme of a subsystem employs updated state measurements of itself and its neighbors to solve the corresponding COCP at each control step. Throughout this work, let $x_{t}^{i}$ denote the measured state of subsystem $i$ at time $t$. We denote by $x_{k,t}^{i}$ the state prediction of subsystem $i$ at time $t+k$, with initial condition $x_{t}^i$. In particular, $x_{0,t}^i\triangleq x_t^i$. $u_{k,t}^i$ follows the same notation for inputs. Let a bold charecter denote the collection of that variable over the prediction horizon, e.g. $\bm x_t^i=[x_{0,t}^i{}^T,\dots,x_{N,t}^i{}^T]^T$, $\bm u_{t}^i=[u_{0,t}^i{}^T,\dots,u_{N-1,t}^i{}^T]^T$. Let $\Psi\triangleq [x^1{}^T,\dots,x^{M}{}^T]^T$. For notation simplicity, the superscript $i$ denoting the subsystem index will be eliminated when only the dynamic of a single subsystem is discussed.

Each subsystem solves a local finite-horizon optimization problem to achieve cooperation. $\forall i=1,\dots, M$, the DMPC scheme of subsytem $i$ is given as follows.
\begin{subequations}
\label{eqn:costi}
\begin{align}
 J_{i}^*(X_{0,t}^i)&= \min_{\bm u_{t}^{i}}\sum_{k=0}^{N-1}l^i_k(X_{k,t}^i, u_{k,t}^i)+l_N^i( X_{N,t}^i)\\
s.t.\quad&  x_{k+1,t}^{i}=A x_{k,t}^{i}+B u_{k,t}^{i}\\
& x_{k,t}^{i}\in\mathcal X^i,\quad  x_{N,t}^{i}\in \mathcal X^i_f, \quad u_{k,t}^{i}\in \mathcal U^i\\
& x_{k+1,t}^{j}=A x_{k,t}^{j}+B u_{k,t}^{j} \quad \forall j\in \mathcal N^i\\
& x_{k,t}^{j}\in\mathcal X^j,\quad  x_{N,t}^{j}\in \mathcal X^j_f, \quad u_{k,t}^{j}\in \mathcal U^j\\
&\forall t\in \mathbb Z^+,\quad i\in \{1,\cdots,M\}
\end{align}
\end{subequations}
where $N\geq 1$ is the prediction horizon, $l_k^i$ and $l_N^i$ are the process cost and the final penalty, respectively. \textcolor{temp}{Define $x_{k,t}^{-i} = [x_{k,t}^{j_1}{}^T, \dots,  x_{k,t}^{j_O}{}^T]^T$ and $X_{k,t}^i=[x_{k,t}^{i}{}^T,  x_{k,t}^{-i}{}^T]^T$ where $O=|\mathcal N^i|$, and let $J_{\Sigma}=\sum_{i=1}^M J_i$.} $\mathcal X_f^i\subseteq \mathcal X^i$ is the terminal constraint set of subsystem $i$.
At time $t$, each subsystem $i,\forall i=1,\dots,M$ solves (\ref{eqn:costi}) for $\bm{u}_t^{i*}$. The first term of the optimized input sequence $\bm{u}_t^{i*}$ is applied to subsystem $i$ with other predicted inputs discarded. Then at next step, all subsystems solve the COCP (\ref{eqn:costi}) again with updated states and updated trajectory assumptions of neighbors. \textcolor{temp}{Without lose of generality, the cost function (\ref{eqn:costi}a) will consider particularly the state difference $||\zeta^i-\zeta^j-d_{i,j}||=||x^i-x^j||$, where $d_{i,j} = \hat{\zeta}^i-\hat{\zeta}^j$, $\forall i=1,\dots, M, \forall j\in\mathcal N^i$. Minimizing the state difference term helps maintaining the rigid grid property of the overall system, which is beneficial for synchronization of subsystems. Details of formulating $l_k(\cdot)$ and $l_N(\cdot)$ taking into account $||x^i-x^j||$ will be given in section \ref{sec:main_result}.}

\section{Main Results}
\label{sec:main_result}
\subsection{DMPC with Switched Cost Functions}
\textcolor{temp}{As all subsystems solve (\ref{eqn:costi}) in a parallel manner, it is difficult to directly combine $||x_{k,t}^i-x_{k,t}^j||, k\geq 1$ in the cost function. In general, the predictions $\{x_{k,t}^j,k\geq 1\}$ are unavailable to subsystem $i\in\mathcal N^j$ until step $t+1$. A commonly adopted approach is introducing a presumed trajectory $\hat {\bm x}_t^{j}$ (together with assumed inputs $\hat{\bm{u}}_t^j$) of subsystem $j$ for its neighbors to solve their COCP at time $t$ \cite{dunbar2006distributed}. By approximating $||\bm{x}_{t}^i-\bm{x}_{t}^j||$ using $||\bm x_{t}^i-\hat{\bm{x}}_{t}^j||$, $J_{i}(X_{0,t}^{i};\bm{x}_t^{j*})$ is estimated by $J_{i}(X_{0,t}^{i};\hat{\bm{x}}_t^{j})$.} The process cost $l_k^i$, $k=0,\dots,N-1$ and final penalty $l_N^i$ at time $t$ are designed as follows
\begin{subequations}
\label{eqn:costk}
\begin{eqnarray}
l_k^i&=&|| x_{k,t}^{i}||_Q^2+|| u_{k,t}^{i}||_R^2\\
&&+\sum_{j\in \mathcal N^i}|| x_{k,t}^{i}- \hat{x}_{k,t}^{j}||_{Q_e}^2 \\
l_N^i&=&||x_{N,t}^{i}||_P^2+\sum_{j\in \mathcal N^i}||x_{N,t}^{i}-\hat x_{N,t}^{j}||_{P_e}^2
\end{eqnarray}
\end{subequations}
where $\hat x_{k+1,t}^j=A\hat x_{k,t}^j+B\hat u_{k,t}^j, k=0,\dots, N-1$. $P$, $P_e$, $Q$, $Q_e$, and $R$ are matrices of appropriate dimensions. \textcolor{temp}{In particular, (\ref{eqn:costk}a) represents the state and input costs of subsystem $i$; (\ref{eqn:costk}b) represents the penalty on state differences between subsystem $i$ and its neighbors $j\in\mathcal N^i$. (\ref{eqn:costk}c) follows the same formulation for the terminal state and terminal state difference.}

For non-iterative DMPC schemes of the form (\ref{eqn:costi})+(\ref{eqn:costk}) with parallel information exchange, $\hat{\bm{x}}_t^i, \forall i=1,\dots, M$ is typically selected based on the optimized state trajectory of the previous step, i.e. $\hat{\bm{x}}_t^i= [x_{1,t-1}^{i*}{}^T,\dots,x_{N,t-1}^{i*}{}^T,\hat{x}_{N,t}^i{}^T]^T$ with $\hat{x}_{N,t}^i= Ax_{N,t-1}^{i*}+B\hat u_{N-1,t}^i$. $\hat u_{N-1,t}$ is obtained using a terminal stabilizing law as in MPC \cite{borrelli2012predictive}. Then, a compatibility constraint on the mismatch between the presumed trajectory and the predicted trajectory of each subsystem $i=1,\dots, M$ is adopted in the following form.
\begin{equation}
    \label{eqn:compatibility}
    ||x_{k,t}^i-\hat{x}_{k,t}^i|| \leq c(t) \qquad k=0,\dots,N,\ \forall t\in\mathbb Z^+
\end{equation}
where $c(t)$ is the compatibility bound to be determined.

When the condition is satisfied that the optimal cost $J_{\Sigma}^*$ at time $t$ is uniquely determined by the state measurement $\Psi(t)$, the compatibility bound can be selected such that $J_{\Sigma}^*$ is monotonically decreasing along the closed-loop trajectory \cite{keviczky2006decentralized,dunbar2006distributed,li2016distributed}. \textcolor{temp}{Stability of the equilibrium is achieved if the implicit control law $\{u_1^*(\cdot)\}$ is further assumed to be continuous in $\Psi$ at the origin \cite{dunbar2012distributed}. As the implicit control law is obtained from solving the COCP (\ref{eqn:costi})+(\ref{eqn:costk})+(\ref{eqn:compatibility}), the optimized input $\bm u_t^{i*}$ could be a function of both $X_{0,t}^i$ and $\hat{\bm x}_t^{-i}$. Therefore, it is desirable to investigate the properties of the non-iterative DMPC with parallel information exchange when the presumed trajectory $\hat{\bm{x}}_{t}^{-i}$ required by subsystem $i$ is not uniquely defined by $X_{0,t}^i$ at time $t$.}

\textcolor{temp}{One approach to addressing dependence of $J_i^{*}$ on $\hat{\bm x}_t^i$ is constructing $\hat{\bm x}_t^i$ to be a trajectory that is uniquely determined by $x_t^i$. However, such presumed trajectories in general cannot preserve monotonic decreasing of $J_i^{*}$ due to lack of information of the previous step. Adding additional constraints guaranteeing monotonic decreasing to (\ref{eqn:compatibility}) could considerably increase conservatism of the non-iterative DMPC scheme.}
Note that the terminal control law is capable of stabilizing the closed-loop system over the terminal set by design, we adopt a switched cost function together with a control invariant set $\bar{\mathcal X}_0\subset\mathcal X$ to address convergence and stability taking into account the general case where $J_i^{*}$ is determined by both $X_{0,t}^i$ and $\hat{\bm x}_t^{-i}$. The rationale is that the cost function of (\ref{eqn:costi}) switches to the cost function of generic decoupled MPC if the state is inside $\bar{\mathcal X}_0$. Thus, stability of each subsystem in $\bar{\mathcal X}_0$ is ensured by conditions of a generic MPC scheme. Convergence of states out of $\bar{\mathcal X}_0$ is ensured by verifying monotonic decreasing of $J^{*}$ as a function uniquely defined by the augmented state $[x_t^i{}^T,\hat{\bm x}_t^{-i}{}^T]^T$. \textcolor{temp}{Therefore, attractiveness of states out of $\bar{\mathcal X}_0$ is addressed from a trajectory optimization perspective instead of constructing a control Lyapunov function based on $J^*(x)$ and $J_{\Sigma}^*(\Psi)$.}

The decoupled MPC problem that only takes into account the state and input of subsystem $i \in \{1,\dots,M\}$ is given as follows
\begin{subequations}
\label{eqn:mpc-ind}
\begin{align}
J_{i}^{\prime*}(x_t^i)& = \min_{\bm u_{t}^{i}}\sum_{k=0}^{N-1}(||x_{k,t}^i||_Q^2+||u_{k,t}^i||_R^2)+||x_{N,t}^i||_P^2\\
s.t.\quad & x_{k+1,t}^{i} = Ax_{k,t}^{i}+Bu_{k,t}^{i}, \quad \forall t\in\mathbb Z^+\\
& x_{k,t}^{i} \in \bar{\mathcal{X}}_0^i,\quad x_{N,t}^{i,i}\in \mathcal{X}_{f}^i, \quad u_{k,t}^{i} \in \mathcal U^i
\end{align}
\end{subequations}

\begin{assumption}
For every subsystem $i\in\{1,\dots, M\}$:
\begin{enumerate}[(a)]
\item $\mathcal X^i$, $\bar{\mathcal X}_0^i$, $\mathcal X_f^i$, and $\mathcal U^i$ are convex and compact polyhedra containing the origin as an interior point.
\item $Q=Q^T\succ 0$, $R=R^T\succ 0$, $Q_e=q_eI, q_e>0$, There exists a linear state feedback law $K(\cdot)$ and $P\succ 0$ such that
\begin{align}
\label{eqn:lqr}
& -||x||_{P}^2+||x||_Q^2+||u||_R^2+||(Ax+Bu)||_{P}^2\leq 0\notag\\
& s.t.\ Ax+Bu\in \mathcal{X}_f, u=K(x)\in\mathcal U,\forall x\in \mathcal X_f \notag
\end{align}
\item $\mathcal X^i$ (respectively, $\mathcal U^i$) is decoupled with $\mathcal X^j$ (respectively, $\mathcal U^j$), $\forall j\in \{1,\dots, M\}\setminus i$.
\end{enumerate}
\end{assumption}

Assumption 1-(a) indicates that state and input constraints can be formulated into linear inequalities in the numerical optimization process. Assumption 1-(b) follows common assumptions adopted in MPC literature \cite{borrelli2012predictive}. Assumption 1-(c) implies that the state and input constraints of a subsystem are not affected by that of its neighbors. For spatially interconnected systems, such assumption could considerably simplify analysis and synthesis of DMPC. Similar assumptions are commonly made in distributed mobile agent systems \cite{keviczky2008decentralized,borrelli2012predictive,dunbar2012distributed}. For DMPC schemes with coupled state and input constraints existing, please refer to \cite{richards2004decentralized,stewart2010cooperative,conte2012distributed}. \textcolor{temp}{Given $Q_e = q_eI$, let $P_e=P_e^T\succ 0$ be such that
\begin{equation}
\label{eqn:pe}
-||x||_{P_e}^2+||x||_{Q_e}^2+||(A+BK)x||_{P_e}^2<0, \forall x\in\mathbb R^n\backslash \{\bm 0\}.
\end{equation}
One approach searching such $P_e$ is solving a discrete-time Ricatti equation by introducing a nominal $R_e\succ 0$.}



\begin{remark}
$\mathcal{X}_f$ is a convex and compact polyhedron. Denote by $\mathcal{X}_0\triangleq Pre_N(\mathcal X_f)$ the $N$-step backward reachable set of $\mathcal{X}_f$ in $\mathcal X$. Then $\mathcal{X}_0$ is convex, compact, and control invariant.
\end{remark}

After combining the DMPC (\ref{eqn:costi}) with decoupled MPC (\ref{eqn:mpc-ind}), the DMPC scheme with switched cost functions is proposed by setting $\bar{\mathcal{X}}_{0}\supseteq \mathcal X_f$ as a control invariant polyhedron set in $\mathcal{X}_0$. $x_{0}^i\in\mathcal X_0^i$ of subsystem $i$ is first regulated to $\bar{\mathcal X}_0^i)$ by (\ref{eqn:costi}), and then, asymptotic stability of the overall system is achieved by (\ref{eqn:mpc-ind}) after all states $\{x^i,i=1,\dots,M\}$ are attracted in $\prod_{i=1}^M\bar{\mathcal X}_0^i$. $\bar{\mathcal X}_0$ serves as a tuning parameter that contributes to spatial properties. A detailed discussion will be given in section III-B. Let $0^-$ denote the time of the initialization step before $t=0$,  the DMPC scheme with switched cost functions is described as follows

\begin{alg}[DMPC with switched cost functions]
\label{alg:dmpc}
\begin{enumerate}[1)]
\item (\textit{initialization}) At time $0^-$, each subsystem $i= 1,\dots,M$ solves (\ref{eqn:costi})+(\ref{eqn:costk}) by setting $q_e = 0$ and $P_e = \bm 0$, and implements the computed control law $u_{0^-}^i=\bm u_{0,0^-}^{i*}$
\item $\forall t\in\mathbb Z^+$, each subsystem queries the state of subsystem $j\in\mathcal N^i$, and updates the state estimations by setting
\begin{equation}
\label{eqn:state_update}
\begin{aligned}
x_{0,t}^{i} &=x_{t}^i\\
\hat{\bm{x}}_{t}^{j} &=[x_{1:N,t-1}^{j*}{}^T, (A+BK)x_{N,t-1}^{j*}{}^T]^T, \quad\forall j\in\mathcal N^i
\end{aligned}
\end{equation}
\item If $\exists i\in \{1,\dots,M\}$ such that $x_{t}^i\notin \bar{\mathcal{X}}_0^i$, then each subsystem $i$ solves problem (\ref{eqn:costi}) with (\ref{eqn:costk}) using updated states from (\ref{eqn:state_update}); else
    all subsystems implement the decoupled MPC law by solving (\ref{eqn:mpc-ind}), $\forall i=1,\dots,M$.
    \label{enu:1}
\item Each subsystem $i$ implements the first term of the input sequence $\bm u_{t}^{i*}$ computed in \textit{Step} \ref{enu:1}.
    \begin{equation}
    \label{eqn:dmpc-law}
    u_t^i=u_{0,t}^{i*}.
    \end{equation}
\item Each subsystem repeats steps \textit{2-5} at time $t+1$ based on the updated states $\{x_{t+1}^i\}$.
\end{enumerate}
\end{alg}

\begin{lemma}[Recursive feasibility]
\label{lemma:feasibility}
Given the conditions in Assumption 1, $\forall \Psi_0\in \prod_{i=1}^M\mathcal X_0^i$, the DMPC problem applying the control law in \textit{Algorithm 1} is recursively feasible.
\end{lemma}

\begin{IEEEproof}
$\mathcal X_f$ is forward invariant with respect to the closed-loop system $x(t+1)=(A+BK)x(t)$. By definition, $\mathcal X_0$ is the $N$-step backward reachable set of $\mathcal X_f$ in $\mathcal X$. Therefore, $\forall x\in \mathcal X_0$, there exists a control sequence $\{u_{0,0^-},\dots,u_{N-1,0^-}\}$, $u_{k,0^-}\in \mathcal U$ such that $x_N\in \mathcal X_f$. Thus, the initialization step is feasible. Feasibility at $t, \forall t\in\mathbb Z^+$ follows the fact that $\{u_{1,t-1},\dots, u_{N-1,t-1}, (A+BK)x_{N,t-1}\}$ is a feasible solution. \end{IEEEproof}


\begin{lemma}
\label{lemma:decrease}
\textcolor{temp}{Suppose Assumption 1 holds. If for every subsystem $i=1,\dots,M$, $\forall t\in\mathbb Z^+$ and $k=0,\dots,N$,
\begin{subequations}
\label{eqn:decrease}
\begin{align}
& DJ_i(x_t^i,\hat{\bm x}_t^{-i},\bm{u_t^{i*}})\leq -\alpha_i(||x_t^i||)\\
& ||\Delta_{k,t}^{i}||\leq \gamma \min_{j\in\mathcal N^i}||x_{t-1}^i-x_{t-1}^j||^2
\end{align}
\end{subequations}
where $||\Delta_{k,t}^i||\triangleq ||x_{k,t}^i-\hat{x}_{k,t}^i||$,  $\alpha_i(\cdot)$ is a class $\mathcal K$ function, $\gamma\in[0,\infty)$ is a constant. Then, by applying the control law in Algorithm 1, i) $\exists T\in\mathbb Z^+$ such that the overall system converges to $\prod_{i=1}^M \bar{\mathcal X}_0^i$ and stays in $\prod_{i=1}^M \bar{\mathcal X}_0^i$ for all $t\geq T$; ii) $\lim_{t\to\infty}J_i(x_t^i,\hat{\bm x}_t^{-i},\bm{u_t^{i*}})=0, \forall i=1,\dots,M$.}
\end{lemma}
\begin{IEEEproof}
Recursive feasibility of the scheme in Algorithm 1 follows the results in Lemma \ref{lemma:feasibility}. Given $x_0^i\in\mathcal X_0^i, \forall i=1,\dots,M$, $J_i(x_t^i,\hat{\bm{x}}_t^{-i},\bm u_{t}^{i*})$ is finite, nonincreasing, and bounded from below by 0. Therefore, both $J_i(x_t^i,\hat{\bm{x}}_t^{-i},\bm u_{t}^{i*})$ and $J_i(x_{t+1}^i,\hat{\bm{x}}_{t+1}^{-i},\bm u_{t+1}^{i*})$ converge to $\hat J_i\geq 0$ as $t\to\infty$.  $||x_t^i|| = \alpha_i^{-1}(\alpha_i(||x_t^i{}||))\to \bm 0$. Given arbitrary $\varepsilon>0$, there exists a $T>0$ such that $\forall t\geq T-1, \forall i=1,\dots, M: \mathrm{(i)} ||x_t^i-x_t^j||^2\leq\varepsilon/(2N(\gamma+1)),\forall j\in\mathcal N^i; \mathrm{(ii)} ||x_t^i||\leq \varepsilon/2$. For every $k=0,\dots,N$, $||\hat{x}_{k,t}^i||=||\hat{x}_{k,t}^i+x_{t+k}^i-x_{t+k}^i||\leq||x_{t+k}^i-\hat{x}_{k,t}^i||+||x_{t+k}^i||\leq N\frac{\gamma\varepsilon}{2N(\gamma+1)}+\frac{\varepsilon}{2}<\varepsilon.$ Therefore, $\hat{\bm{x}}_t^i\to \bm 0$ as $t\to\infty$, which indicates ${\bm{x}}_t^{i*}\to \bm 0$ as $t\to\infty$. Furthermore, $\lim_{t\to\infty}x_{k+1,t}^{i*} = \lim_{t\to\infty}Ax_{k,t}^{i*} + \lim_{t\to\infty}Bu_{k,t}^{i*}$, and $B$ is of full rank, we have $\mathrm{ker}(B)=\{\bm 0\}$ and $\bm{u}_{t}^{i*}\to 0$ as $t\to\infty$. Thus, $J_i(x_t^i,\hat{\bm{x}}_t^{-i},\bm{u}_t^{i*})\to\hat{J}_i=0$ as $t\to\infty$. Since $J_i(\cdot)$ is a quadratic function of the augmented vector $\chi=[x_t^i{}^T,\hat{\bm{x}}_t^{-i}{}^T,\bm{u}_{t}^{i*}{}^T]^T$, there exists a class $\mathcal K_{\infty}$ function $\underline{\alpha}(||\chi(t)||)$ as a lower bound. Denote by $\partial \bar{\mathcal X}_0^i$ the boundary of $\bar{\mathcal X}_0^i$. Since $\bar{\mathcal X}_0^i$ is compact, there exists $d_{\mathrm{min}}=\min||x||$, $\forall x\in\partial \bar{\mathcal X}_0^i$. Let $\Omega_d$ be the level set of $J_i(\chi)$, $\Omega_d = \{\chi: J_i(\chi)\leq \underline{\alpha}(d_{\mathrm{min}})\}\subseteq\mathcal B_{d_{\min}}(\bm 0)$. Since $\chi$ converges to $\Omega_d$ in finite time, hence $\Phi(t,x_0^i)$ converges to $\bar{\mathcal X}_0^i$ in finite time.
\end{IEEEproof}

\begin{remark}
\textcolor{temp}{Note that the proof of Lemma \ref{lemma:decrease} does not require the cost function uniquely defined by $[x_t^i{}^T,x_t^{-i}{}^T]^T$. The control law can be based on information of both $[x_t^i{}^T,x_t^{-i}{}^T]^T$ and $\hat{\bm x}_t^{-i}$. The property that $\bm x$ converges to $\prod_{i=1}^M\bar{\mathcal X}_0^i$ in finite time is used to enable the switch indicator for all subsystems. In particular, each subsystem saves the same $d_{\mathrm{min}}$. At each control step, $J_i^*$ is compared with $d_{\mathrm{min}}$, $\forall i=1,\dots,M$. A subsystem will send out a \textit{switch-ready} message if $J^*\leq \underline{\alpha}(d_{\mathrm{min}})$. It also records all \textit{switch-ready} messages for all reachable subsystems. A subsystem will then send out a \textit{list-complete} message if it receives switch-ready status of all subsystems. The overall system switches to the decoupled control law if each subsystem receives a list-complete indication confirmed by all other subsystems. The decision consensus process is ensured by the assumption that $\mathcal G$ has a spanning tree. With neighbor-to-neighbor communication constraints, the consensus process may occupy multiple information exchange cycles. In cases where temporal consensus channels are available, consensus on the switching action can be conducted in one information exchange cycle.}
\end{remark}


\begin{lemma}
\textcolor{temp}{Given conditions in Assumption 1, if the state difference for each subsystem $i=1,\dots,M$ satisfies
\begin{subequations}
\label{eqn:condi_last}
\begin{align}
x_{N,t}^{i} =& \hat x_{N,t}^{i}\\
||\Delta_{k,t}^i||\leq & \min_{j\in\mathcal N^i}\frac{||x_{t-1}^i-x_{t-1}^j||^2}{4(N-1)\rho_{\max}},\ \forall k=0,\dots, N
\end{align}
\end{subequations}
}
where $\rho_{\max}$ is an upper bound on $||x^i||, \forall x^i\in\mathcal X^i$ for every $i=1,\dots,M$. Then, $\forall \Psi_0\in\prod_{i=1}^M\mathcal X_0^i$ the overall system converges to $\prod_{i=1}^M\bar{\mathcal X}_0^i$ in finite time $T<\infty$ and stays in $\prod_{i=1}^{M}\bar{\mathcal X}_0^i$ for all $t\geq T$ by applying the control law proposed in \textit{Algorithm 1}.
\end{lemma}

\begin{IEEEproof}
To prove this property, we are going to show that solving (\ref{eqn:costi})+(\ref{eqn:costk}) with condition (\ref{eqn:condi_last}) satisfies conditions (\ref{eqn:decrease}a) and (\ref{eqn:decrease}b) given in Lemma \ref{lemma:decrease}. Since $\mathcal X$ is compact, there exist $\rho_{\mathrm{max}}\triangleq \max_{i=1:M}\max_{x^i\in\mathcal X^i} ||x^i||$ as an upper bound, and $\rho_{\mathrm{max}}<\infty$. Hence, (\ref{eqn:decrease}b) is satisfied. Furthermore, $\hat{\bm x}_t^i$ (together with $\hat{\bm{u}}_t^i$) is a feasible solution to constraints in (\ref{eqn:condi_last}), recursive feasibility of Algorithm 1 holds following results in Lemma 1. For the condition (\ref{eqn:decrease}a), (\ref{eqn:condi_last}a) together with (\ref{eqn:pe}) indicates that $\sum_{j\in\mathcal N^i}||x_{N,t-1}^{i*}-x_{N,t-1}^{j*}||_{Q_e}^2+\sum_{j\in\mathcal N^i}||(A+BK_l)(x_{N,t-1}^{i*}-x_{N,t-1}^{j*})||_{P_e}^2-\sum_{j\in\mathcal N^i}||x_{N,t-1}^{i*}-x_{N,t-1}^{j*}||_{P_e}^2\leq 0$.
\begin{subequations}
\label{eqn:dj}
\begin{align}
& J_i(x_{k,t+1}^{i},\hat{\bm{x}}_{t+1}^{-i},\bm{u}_{t+1}^{i*})-J_i(x_{k,t}^{i},\hat{\bm{x}}_{t}^{-i},\bm{u}_{t}^{i*})\\
\leq & J_i(x_{k,t+1}^{i},\hat{\bm{x}}_{t+1}^{-i},\hat{\bm{u}}_{t+1}^{i})-J_i(x_{k,t}^{i},\hat{\bm{x}}_{t}^{-i},\bm{u}_{t}^{i*})\\
\leq & -\sum_{j\in\mathcal N^i}||x_{t}^{i}-x_{t}^{j}||_{Q_e}^2 - ||x_{t}^i||_{Q}^2 - ||u_{0,t}^{i*}||_{R}^2\\
& + \sum_{k=1}^{N-1}\sum_{j\in\mathcal N^i}(||x_{k,t}^{i*}-x_{k,t}^{j*}||_{Q_e}^2-||\Delta_{k,t}^{i}||_{Q_e}^2)
\end{align}
\end{subequations}

(\ref{eqn:dj}d) is further expanded as follows
\begin{align}
(\ref{eqn:dj}\mathrm{d})=&\sum_{k=1}^{N-1}\sum_{j\in\mathcal N^i}\left(q_e(||x_{k,t}^{i*}-x_{k,t}^{j*}||+||x_{k,t}^{i*}-\hat x_{k,t}^{j}||)\right.\notag\\
&\qquad\left.(||x_{k,t}^{i*}-x_{k,t}^{j*}||-||x_{k,t}^{i*}-\hat x_{k,t}^{j}||)\right)\notag\\
\leq & \sum_{k=1}^{N-1}\sum_{j\in\mathcal N^i}\left(q_e(||x_{k,t}^{i*}-x_{k,t}^{j*}||+||x_{k,t}^{i*}-\hat x_{k,t}^{j}||)(||\Delta_{k,t}^j||)\right)\notag\\
\leq & \sum_{k=1}^{N-1}\sum_{j\in\mathcal N^i}4q_e\rho_{\max}||\Delta_{k,t}^j||\notag
\end{align}

Further implementing (\ref{eqn:condi_last}b) to (\ref{eqn:dj}d) indicates $(\ref{eqn:dj}\mathrm{a})\leq - ||x_{t}^i||_{Q}^2 - ||u_{0,t}^{i*}||_{R}^2$, which meets the condition of (\ref{eqn:decrease}a). Therefore, the overall system converges to $\prod_{i=1}^M\bar{\mathcal X}_0^i$ in finite time.
\end{IEEEproof}

Combined with the compatibility constraints proposed in (\ref{eqn:condi_last}), the DMPC scheme proposed in Algorithm 1 is further described as follows.


\begin{alg}
\label{alg:dmpc2}
\begin{enumerate}[1)]
\item Follow steps \textit{1-2} in Algorithm 1
\item If $\exists i\in\{1,\dots,M\}$ such that $x_t^i\notin \mathcal X_f$, then each subsystem solves (\ref{eqn:costi})+(\ref{eqn:costk})+(\ref{eqn:condi_last}) using updated states from (\ref{eqn:state_update}); else all subsystem implement the decoupled MPC law. \label{enu:2b}
\item Follow steps \textit{4-5} in Algorithm 1
\end{enumerate}
\end{alg}

\begin{lemma}[{\cite[page~133]{borrelli2012predictive}}]
Consider the mp-QP (\ref{eqn:mpqp}) and let $H\succ 0$. Then the optimizer $z^*(x):\mathcal K\to \mathbb R^s$ is continuous and piecewise affine on polyhedra, where $\mathcal K$ is the feasible set. In particular it is affine in each critical region.
\end{lemma}

\begin{theorem}
Given the conditions in Assumption 1 and Lemma 3. Then, by applying the DMPC law in Algorithm 2
\begin{enumerate}[(i)]
\item The origin of the overall system comprising subsystems (\ref{eqn:dynamic})  is a Lyapunov stable equilibrium; $\forall \Psi_0\in \prod \mathcal X_0^i$, $i=1,\dots,M$, $\lim_{t\to\infty}\Phi(t,\Psi_0)=\bm 0$.
\item The DMPC scheme in Algorithm 2 is a series of quadratic programs that can be transformed into the following form
    \begin{equation}
    \label{eqn:cocp}
    \begin{aligned}
    J_i^*(U_0^i, \bm u_t^i)=&\min_{\bm u_t^i}\begin{bmatrix}
    \bm u_t^i{}^T & U_0^i{}^T
    \end{bmatrix}
    \begin{bmatrix}
    H & F^T \\
    F & Y
    \end{bmatrix}
    \begin{bmatrix}
    \bm u_t^i{}^T & U_0^i{}^T
    \end{bmatrix}^T\\
    & s.t. \quad G\bm u_t^i\leq w + EU_0^i
    \end{aligned}
    \end{equation}
    where $\bm u_t^i$ is the optimization variable, $H$, $F$, and $Y$ are matrix coefficients uniquely determined by the augmented initial state $U_0^i\triangleq[x_t^i{}^T,\bm{x}_t^{-i}{}^T]^T$. The optimizer $\bm{u}_t^{i*}$ is a piecewise affine function of the augmented initial condition $Z_0^i$.
\end{enumerate}
\end{theorem}

\begin{IEEEproof}
(i) Given the conditions in Assumption 1 and $\bar{\mathcal X}_0$, the decoupled MPC schemes respectively stabilize each subsystem after cost function switching, which ensures Lyapunov stability of the overall system. The decoupled MPC schemes further ensures convergence for every $\Psi\in\prod_{i=1}^M \bar{\mathcal X}_0^i$. Convergence of $\Phi(t,\Psi_0), \forall \Psi_0\in\prod_{i=1}^M\mathcal X_0$ to $\prod_{i=1}^M\bar{\mathcal X}_0^i$ follows the results in Lemma 3.


(ii) For the first stage when $x_0^i\notin \bar{\mathcal X}_0^i$, $\bm x_t^i$ is affine of $\bm u_t^i$, and both the state and input constraints are convex polyhedra. Moreover, (\ref{eqn:condi_last}b) can be ensured by the following inequality
\begin{equation}
    ||\Delta_{k,t}^i||_{\infty}\leq \min_{j\in\mathcal N^i}\frac{||x_{t-1}^i-x_{t-1}^j||^2}{4\sqrt{n}(N-1)\rho_{\max}},\ \forall k=0,\dots, N
\end{equation}
Therefore, the COCP (\ref{eqn:costi})+(\ref{eqn:costk})+(\ref{eqn:condi_last}) is a quadratic program as described by (\ref{eqn:cocp}). By defining $z=\bm u_t+H^{-1}F^TU_0^i$, $z\in \mathbb R^s$, and removing $U_0^i{}^TYU_0^i$ that is independent of $\bm u_t^{i*}$, (\ref{eqn:cocp}) is transformed into the following form
\begin{equation}
\begin{aligned}
\hat{J}^*(X_0^i) &= \min_{z}z^THz\\
& s.t.\quad Gz\leq w+ SU_0^i
\end{aligned}
\end{equation}
where $S=E+GH^{-1}F^T$. Following Lemma 4, the optimizer $z^*(U_0^i)$ is continuous and piecewise affine on polyhedra.

For the second stage, the DMPC scheme changes into $M$ decoupled MPC problems. The problem is quadratic program following the analysis of generic MPC schemes (see, e.g., \cite{borrelli2012predictive}).
\end{IEEEproof}

\begin{remark}
\textcolor{temp}{The non-iterative DMPC scheme in Algorithm 2 with constraint (\ref{eqn:condi_last}) is different from a Lyapunov-based approach. In the proposed cost-switching approach, we only require convergence of the DMPC scheme (\ref{eqn:costi})+(\ref{eqn:costk})+(\ref{eqn:condi_last}). No constraints are required to ensure that $J_{\Sigma}^*(\cdot)$ is a control Lyapunov function. Therefore, the COCP is addressed from a trajectory optimization perspective that implements results of the previous optimization step. The DMPC control law is generated based on both the current state of a subsystem $i$ and the state assumptions of its neighbors $j\in\mathcal N^i$, $\forall i=1,\dots, M$.}
\end{remark}

\subsection{Independent Obstacle-Avoidance Capability}
The implementation of cost-switching strategy in the DMPC scheme introduces decoupled regulation capability of subsystems in $\bar{\mathcal X}_0$. Therefore, it is natural to investigate independent obstacle-avoidance capability of the decoupled subsystems. In this section, let $\tilde{\mathcal B}_r(\zeta)$ denote the ball of radius $r$ centered in $\zeta$ in the absolute coordinates. Let $\perp_{s}$ denote the projection operator of a set on the spatial subspace. In particular, we have the following result.

\begin{figure}
\centering
\includegraphics[width=5.8cm]{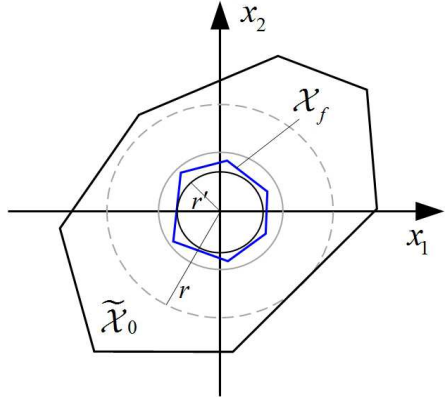}
\caption{A 2-dimonsional diagram of terminal sets of the DMPC scheme.}
\label{fig:sets}
\end{figure}

\begin{theorem}
There exists a configuration of $\bar{\mathcal{X}}_0\subseteq\mathcal X_0$ for each subsystem, such that i) \textcolor{temp}{the subsystems applying control law in Algorithm 2 with initial condition $x_0\in\bar{\mathcal{X}}_0$ have no intersections in the spatial subspace, i.e. $\cap_{i=1}^M \mathcal S_0^i=\emptyset$. ii) $\exists r\prime >0$ such that the state shift from $\bm 0$ to an avoidance state $x^i\in\mathcal B_{r^{\prime}}(\bm 0)$ and back to $\bm 0$ by change of coordinates and using (\ref{eqn:mpc-ind}) stays in $\bar{\mathcal X}_0$.}
\end{theorem}
\begin{IEEEproof}
i) Since the systems are spatially distributed, given the system configuration, there exists a polyhedron set $\tilde{\Xi}_0^{i}, i=1,\dots, M$ such that $\forall j\in\mathcal N^i, \perp_s(\tilde{\Xi}_0^i)\cap\perp_s(\tilde{\Xi}_0^j)=\emptyset$, where $\tilde\Xi_0=\{\zeta:\zeta=x+\hat\zeta, x\in\tilde{\mathcal X}_0\}$. There exists $r>0$ such that $\mathcal B_{r}(\bm 0)\subset\tilde{\mathcal{X}}_0^i$. Let $\partial\mathcal X_f^i$ denote the boundary of $\mathcal X_f^i$. Define $\bar J_i=\max_{x\in\partial\mathcal X_f^i}J_i^{*}(x^i,\hat{\bm{x}}^{-i})$. Let $\Omega_d$ be the level set of $J_i^{*}(x^i,\hat{\bm{x}}^{-i})$, $\Omega_d=\{x:J_i^{*}(x^i,\hat{\bm{x}}^{-i})\leq \bar J_i\}$. Further choose $\mathcal{X}_{f}^{i}$ such that $\Omega_d\subseteq\mathcal B_{\frac{r}{2}}(\bm 0)$. Let $\bar{\mathcal X}_0^i=Pre_N(\mathcal X_f^i)$ in $\tilde{\mathcal X}_0^i$. $\bar{\mathcal X}_0^i$ is a control invariant set with respect to (\ref{eqn:mpc-ind}), and $\forall j\in\mathcal N^j$, $\mathcal{S}_0^{i}\cap \mathcal{S}_0^{j}=\emptyset$. By implementing the DMPC scheme in Algorithm 2, each subsystem is controlled by (\ref{eqn:mpc-ind}) if conditions at step (2) are reached. Therefore, the trajectories have no intersections in the spatial subspace.

ii) Select $\mathcal B_{r\prime}(\bm 0)\subset\mathcal{X}_f, r\prime>0$. Denote by $\mathcal Z_f^i(x^i)$ the terminal set $\mathcal X_f^i$ shifted to $x^i\in\mathcal B_{r\prime}^i(\bm 0)$ by change of coordinates, $\forall i\in\{1,\dots,M\}$. Then for each $x^i\in\mathcal B_{r\prime}^i(\bm 0)\subset\mathcal{X}_f^i$,
\begin{subequations}
\label{eqn:intres}
\begin{align}
&\quad x^i\in\mathcal B_{r\prime}^i(\bm 0) \to x=\bm 0\in\mathcal{Z}_f^i(x^i)\\
&\quad \mathcal{Z}_f^i(x^i)\subset \mathcal B_{r/2}(x^i) \subset\mathcal B_{r}^i(\bm 0), \quad\forall x^i\in\mathcal{B}_{r\prime}^i(\bm 0)
\end{align}
\end{subequations}
(\ref{eqn:intres}a) indicates that state shift form $\bm 0$ to $x^i$ using (\ref{eqn:mpc-ind}) is feasible, and the trajectory will stay in $\mathcal B_{r/2}(x^i)$. (\ref{eqn:intres}b) indicates that, given initial condition $x_0\in\mathcal{B}_{r\prime}(\bm 0)$, the closed-loop trajectory will stay in $\mathcal B_r(\bm 0)$, $\forall t\in\mathbb Z^+$. Therefore, each subsystem $i=1,\dots,M$ applying (\ref{eqn:mpc-ind}) in Algorithm 2 for state shifting will stay in $\bar{\mathcal X}_0^i, \forall t\in\mathbb Z^+$.
\end{IEEEproof}

\begin{remark}
Theorem 2-i shows that the subsystems have no intersections in the spatial subspace after cost function switching to (\ref{eqn:mpc-ind}a). In addition, the distributed subsystems can independently regulate spatial disturbances at this stage. Thus, safety-related spatial constraints only need to be added to (\ref{eqn:condi_last}) for the first stage if collision-free properties are required. For an illustrating example in flexible platooning, please refer to \cite{liu2017non-iterative}. Note that the feasible set accounting for safety constraints may not necessarily be $\prod_{i=1}^M Pre_N(\mathcal X_f^i)$. Theorem 2-ii indicates that each subsystem has the potential of independently avoiding stationary obstacle if the error dynamic of the overall system is at rest. An independent obstacle avoidance can be processed if such avoidance manner can be addressed by state shifting in $\mathcal B_{r\prime}(\bm 0)$.
\end{remark}

\section{Numerical Results}
\begin{figure}
\centering
\includegraphics[width=8cm]{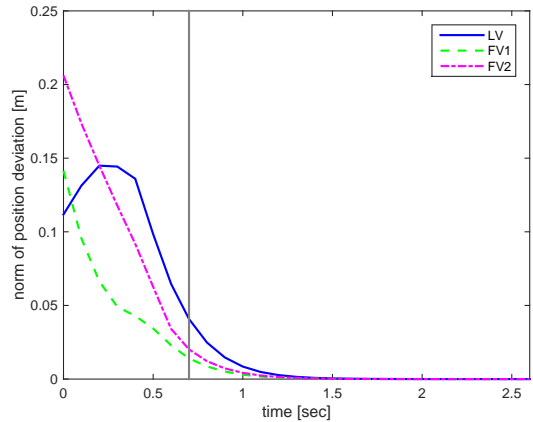}
\caption{Vehicle formation using DMPC with switched cost functions. The vertical gray line indicates the switching time of cost functions.}
\label{fig:tac16eg}
\end{figure}

\begin{figure}
\centering
\includegraphics[width=8cm]{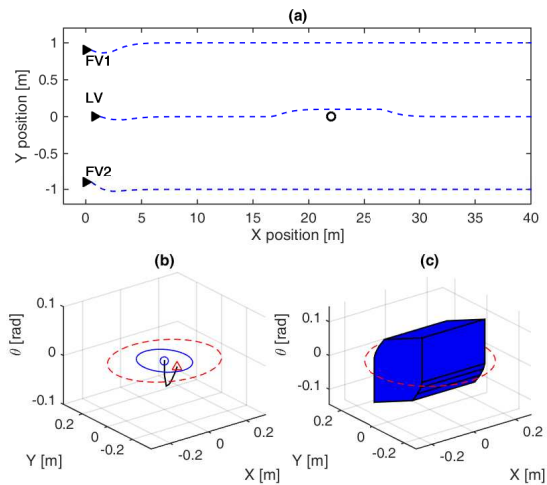}
\caption{(a) vehicle trajectories with the lead vehicle avoiding an stationary obstacle individually; (b) state trajectory (black) to a avoidance state (denoted as a red triangle). The blue circle and the dashed red circle depict the $r^{\prime}$-ball and the $r/2$-ball in $X-Y$ plane, respectively; (c) the terminal set $\mathcal X_f$ with the dashed red circle that is the same as denoted in (b).}
\label{fig:path}
\end{figure}

Consider an unmanned ground vehicle (UGV) formation of three vehicles with the following planar kinematic model
\begin{equation}
\left\{\begin{array}{l}
\dot s=v\cos{\theta}\\
\dot y=v\sin{\theta}\\
\dot\theta=\omega
\end{array}\right.
\end{equation}
where $s$, $y$ and $\theta$ are vehicle longitudinal position, lateral position, and heading angle, respectively. $v$ and $\omega$ are vehicle speed and steering rate, respectively. \textcolor{temp}{The linearized discrete-time model is given as follows by decoupling the longitudinal and lateral dynamics}
\begin{equation}
    \zeta(t+1)=A\zeta(t)+B\varrho(t)
\end{equation}
where
\begin{equation*}
A=\begin{bmatrix}
1 & 0 & 0\\
0 & 1 & 0.5\\
0 & 0 & 1\\
\end{bmatrix}, B=\begin{bmatrix}
0.1 & 0\\
0 & 0\\
0 & 0.1\\
\end{bmatrix}
\end{equation*}
with $\zeta=[s,y,\theta]^T$, $\varrho=[v,\omega]^T$. \textcolor{temp}{Let $\hat\zeta(t)=[s(0)+\bar v t, 0,0]^T$ and $[\bar v, 0]^T$ denote the reference state and the reference input of the reference frame, respectively. $\bar v = 5m/s$ is the cruise speed, $x=\zeta-\hat \zeta$, $u=\varrho-\hat\varrho$. The sampling interval is selected as $0.1s$. The DMPC scheme implements the error dynamic model $x(t+1)=Ax(t)+Bu(t)$ to regulate state deviation from the reference state.} All three vehicles are set with initial state deviations. In the DMPC scheme, cost functions switch to decoupled indexes at $t=0.7s$. The control problem switches from (\ref{eqn:costi}) to (\ref{eqn:mpc-ind}). Fig. \ref{fig:tac16eg} shows norm of state deviations of three vehicles. State deviations converge to zero after 2 seconds. Fig. \ref{fig:path} shows the results regarding independent obstacle-avoidance capability. In this case, the lead vehicle encounters a stationary obstacle that can be avoided by a state transition in the terminal set. Since the control laws are decoupled for states inside $\bar{\mathcal X}_0$, the lead vehicle individually generates an obstacle-avoiding path without disturbing the followers. Additionally, states of all vehicles are inside $\tilde{\mathcal X}_0$, the obstacle avoidance maneuver is collision-free for the formation.



%


\section{Conclusions}
In this note, a non-iterative DMPC scheme for spatially interconnected systems has been investigated implementing cost function switching strategy. With cost-switching, stability of the overall system is guaranteed by the terminal decoupled MPC law that does not require compatibility constraints. Convergence of the closed-loop trajectory and convergence of the optimal cost have been investigated by taking into account the causal link between the presumed trajectory and the optimized trajectory. Compatibility constraints are presented based on convergence analysis of the sequence of the augmented state in the cost function. Independent obstacle-avoidance capability is also tested by taking advantage of the decoupled terminal sets.



\ifCLASSOPTIONcaptionsoff
  \newpage
\fi



\bibliographystyle{IEEEtran}
\bibliography{IEEEabrv,RefListPeng}

\end{document}